# ANALYSIS OF CAPUTO LINEAR FRACTIONAL DYNAMIC SYSTEMS WITH TIME DELAYS THROUGH FIXED POINT THEORY

M. De la Sen

**Faculty of Science and Technology**. University of the Basque Country

Leioa( Bizkaia). Aptdo. 644 de Bilbao. 48080- Bilbao. SPAIN

**Abstract**. This paper investigates the global stability and the global asymptotic stability independent of the sizes of the delays of linear time-varying Caputo fractional dynamic systems of real fractional order possessing internal point delays. The investigation is performed via fixed point theory in a complete metric space by defining appropriate non-expansive or contractive self- mappings from initial conditions to points of the state- trajectory solution. The existence of a unique fixed point leading to a globally asymptotically stable equilibrium point is investigated in particular under easily testable sufficiency-type stability conditions. The study is performed for both the uncontrolled case and the controlled case under a wide class of state feedback laws.

**Keywords**. Caputo fractional systems Time- delay systems, Point delays, Fixed point theory, Stability.

## 1. Introduction

Fractional calculus is concerned with the calculus of integrals and derivatives of any arbitrary real or complex orders. In this sense, it may be considered as a generalization of classical calculus which is included in the theory as a particular case. There is a good compendium of related results with examples and case studies in [1]. Also, there is an existing collection of results in the background literature concerning the exact and approximate solutions of fractional differential equations of Riemann – Liouville and Caputo types, [1-4], fractional derivatives involving products of polynomials, [5-6], fractional derivatives and fractional powers of operators, [7-9], boundary value problems concerning fractional calculus (see, for instance , [1] , [10]) etc . On the other hand, there is also an increasing interest in the recent mathematical related to dynamic fractional differential systems oriented towards several fields of Science like Physics, Chemistry or Control Theory. Perhaps the reason of interest of fractional calculus is that, the numerical value of the fraction parameter allows a closer characterization of eventual uncertainties present in the dynamic model. We can also find, in particular, abundant literature concerned with the development of Lagrangian and Hamiltonian formulations where the motion integrals are calculated though fractional calculus and also in related investigations concerned dynamic and damped and diffusive systems [11-17] as well as the characterization of impulsive responses or its use in Applied Optics related, for instance, to the formalism of fractional derivative Fourier plane filters (see, for instance, [16-18]) and Finance [19]. Fractional calculus is also of interest in Control Theory concerning for instance, heat transfer, lossless transmission lines, the use of discretizing devices supported by fractional calculus, etc. (see, for instance [20-22]). In particular, there are several recent applications of fractional calculus in the fields of filter design, circuit theory and robotics, [21-22], and signal processing, [17]. Fortunately, there is an increasing mathematical literature currently available on fractional differ-integral calculus which can formally support successfully the investigations in other related disciplines.

This paper is concerned with the investigation of the solutions of time-invariant fractional differential dynamic systems, [23-24], involving point delays what leads to a formalism of a class of functional differential equations, [25-31]. Functional equations involving point delays are a crucial mathematical tool to investigate real process where delays appear in a natural way like, for instance, transportation problems, war and peace problems or biological and medical processes. The main interest of this paper is concerned with the positivity and stability of solutions independent of the sizes of the delays and also being independent of eventual coincidence of some values of delays if those ones are, in particular, multiple related to the associate matrices of dynamics. Most of the results are centred in characterizations via Caputo fractional differentiation although some extensions are presented concerned with the classical Riemann- Liouville differ-integration. It is proved that the existence nonnegative solutions independent of the sizes of the delays and the stability properties of linear time-invariant fractional dynamic differential systems subject to point delays may be characterized with sets of precise mathematical results.

On the other hand, fixed point theory is a very powerful mathematical tool to be used in many applications where stability knowledge is needed. For instance, the concepts of contractive, weak contractive, asymptotic contractive and non-expansive mappings have been investigated in detail in many papers from several decades ago (see, for instance, [32-34] and references therein) . It has been found, for instance, that contractivity, weak contractivity and asymptotic contractivity ensure the existence of a unique fixed pointing complete metric or Banach spaces. Some theory and applications of some types of functional equations in the context of fixed point theory have been investigated in [35-36]. Fixed point theory has also been employed successfully in stability problems of dynamic systems such as time- delay and continuous-time/digital hybrid systems and in those involving switches among different parameterizations. This paper is concerned with the investigation of fixed points in Caputo linear fractional dynamic systems of real order $\alpha$ which involved delayed dynamics subject to a finite set of bounded point delays which can be of arbitrary sizes. The self-mapping defined in the state space from initial conditions to points of the state – trajectory solution are characterized either as non-expansive or as contractive. The first case allows to establish global stability results while the second one characterizes global asymptotic stability.

**1.1. Notation**. $\mathbf{C}$, $\mathbf{R}$ and $\mathbf{Z}$ are the sets of complex, real and integer numbers, respectively.

$\mathbf{R}_+$ and $\mathbf{Z}_+$ are the sets of positive real and integer numbers, respectively; and $\mathbf{C}_+$ is the set of complex numbers with positive real part.

$\mathbf{C}_{0+} := \mathbf{C}_+ \cup \{i\omega : \omega \in \mathbf{R}\}$, where $\mathbf{i}$ is the complex unity, $\mathbf{R}_{0+} := \mathbf{R}_+ \cup \{0\}$ and $\mathbf{Z}_{0+} := \mathbf{Z}_+ \cup \{0\}$.

$\mathbf{R}_-$ and $\mathbf{Z}_-$ are the sets of negative real and integer numbers, respectively; and $\mathbf{C}_-$ is the set of complex numbers with negative real part.

$\mathbf{C}_{0-} := \mathbf{C}_- \cup \{i\omega : \omega \in \mathbf{R}\}$, where $\mathbf{i}$ is the complex unity, $\mathbf{R}_{0-} := \mathbf{R}_- \cup \{0\}$ and $\mathbf{Z}_{0-} := \mathbf{Z}_- \cup \{0\}$.

$\overline{N} := \{1, 2, ..., N\} \subset \mathbf{Z}_{0+}$, " $\vee$ "is the logic disjunction, and " $\wedge$ "is the logic conjunction. $[t/h]$ is the integer part of the rational quotient t /h .

$\sigma(M)$ denotes the spectrum of the real or complex square matrix M ( i.e. its set of distinct eigenvalues).

$\| \ \|$ denotes any vector or induced matrix norm. Also, $\|m\|_p$ and $\|M\|_p$ are the $\ell_p$-norms of the vector m or (induced) real or complex matrix M, and $\mu_p(M)$ denote the $\ell_p$ measure of the square matrix M, [20]. The matrix measure $\mu_p(M)$ is defined as the existing limit

$$\mu_p(M) := \lim_{\varepsilon \to 0^+} \frac{\|I_n + \varepsilon X\|_p - \varepsilon}{\varepsilon}$$

which has the property

$$max\left(-\|M\|_p, \max_{i \in \bar{n}} re\, \lambda_i(M)\right) \leq \mu_p(M) \leq \|M\|_p$$

for any square n-matrix M of spectrum $\sigma(M) = \{\lambda_i(M) \in C : 1 \leq i \leq n\}$. An important property for the investigation of this paper is that $\mu_2(M) < 0$ if M is a stability matrix: i.e. if $re\, \lambda_i(M) < 0\,;\, 1 \leq i \leq n$.

$\|\ \|_\infty$ denotes the supremum norm on $R_{0+}$, or its induced supremum metric, for functions or vector and matrix functions without specification of any point-wise particular vector or matrix norm for each $t \in R_{0+}$. If point-wise vector or matrix norms are specified, the corresponding particular supremun norms are defined by using an extra subscript. Thus, $\|m\|_{p\infty} := \sup_{t \in R_{0+}} \|m(t)\|_p$ and $\|M\|_{p\infty} := \sup_{t \in R_{0+}} \|M(t)\|_p$ are, respectively, the supremum norms on $R_{0+}$ for vector and matrix functions of domains in $R_{0+} \times R^n$, respectively, in $R_{0+} \times R^{n \times m}$ defined from their $\ell_p$ point-wise respective norms for each $t \in R_{0+}$.

$I_n$ is the n-th identity matrix.

$K_p(M)$ is the condition number of the matrix M with respect to the $\ell_p$-norm.

$\bar{k} := \{1, 2, \cdots, k\}$

The sets $BPC^{(i)}(dom, codom)$ and $PC^{(i)}(dom, codom)$ are the sets of functions of a certain domain and codomain which are of class $C^{(i-1)}(dom, codom)$ and with the i-th derivative is bounded piecewise continuous, respectively, piecewise continuous in the definition domain.

## 2. Caputo fractional linear dynamic systems with point constant delays and the contraction mapping theorem

Consider the linear functional Caputo fractional dynamic system of order $\alpha$ with r delays:

$$\left(D_{0+}^{\alpha} x\right)(t) := \frac{1}{\Gamma(k-\alpha)} \int_0^t \frac{x^{(k)}(\tau)}{(t-\tau)^{\alpha+1-k}} d\tau$$

$$= \sum_{i=0}^{r} \hat{A}_i(t) x(t-r_i) + B(t) u(t) = \sum_{i=0}^{r} A_i x(t-r_i) + \sum_{i=0}^{r} \tilde{A}_i(t) x(t-r_i) + B(t) u(t)$$

(2.1)

with $k - 1 < \alpha (\in R_+) \leq k\,;\quad k-1, k \in Z_{0+}$, $0 = r_0 < r_1 < r_2 < \ldots < r_r = h < \infty$ being distinct constant delays, where $r_i\,(i \in \bar{r})$ are the r (in general incommensurate delays) $0 = r_0 < r_i\,(i \in \bar{r})$ subject to the system piecewise continuous bounded matrix functions of delayed dynamics $\hat{A}_i : \mathbf{R}_{0+} \to \mathbf{R}^{n \times n}$ ($i \in \bar{r} \cup \{0\}$) which are decomposable as a (non unique) sum of a constant matrix plus a bounded matrix function of time ; i.e. $\hat{A}_i(t) = A_i + \tilde{A}_i(t)$, $\forall t \in \mathbf{R}_{0+}$, and $B : \mathbf{R}_{0+} \to \mathbf{R}^{n \times m}$ is the piecewise continuous bounded control matrix. The initial condition is given by k n-real vector functions $\varphi_j : [-h, 0] \to \mathbf{R}^n$, with $j \in \overline{k-1} \cup \{0\}$, which are absolutely continuous except eventually in a set of zero measure of $[-h, 0] \subset \mathbf{R}$ of bounded discontinuities with $\varphi_j(0) = x_j(0) = x^{(j)}(0) = x_{j0}$, $j \in \overline{k-1} \cup \{0\}$. The function vector $u : \mathbf{R}_{0+} \to \mathbf{R}^m$ is any given bounded piecewise continuous control function. The following result is concerned with the unique solution on $\mathbf{R}_{0+}$ of the above differential fractional system (3.1). The proof, which is based on Picard-Lindeloff theorem, follows directly from a parallel existing result from the background literature on fractional differential systems by grouping all the additive forcing terms of (2.1) in a unique one (see, for instance [1], Eqs. (1.8.17), (3.1.34)-(3.1.49), with $f(t) \equiv \sum_{i=1}^{r} A_i x(t - h_i) + Bu(t)$). For the sake of simplicity, the domains of initial conditions and controls are all extended to $[-h, 0] \cup \mathbf{R}_{0+}$ by zeroing them on the irrelevant intervals of $[-h, 0)$ so that any solution for $t \in \mathbf{R}_{0+}$ of (2.1) is identical to the corresponding one under the above given definition domains of vector functions of initial conditions and controls.

**Theorem 2.1**. The linear and time-invariant differential functional fractional dynamic system (2.1) of any order $\alpha \in C_{0+}$ has a unique continuous solution on $[-h, 0] \cup \mathbf{R}_{0+}$ satisfying:

(a) $x \equiv \varphi \equiv \sum_{j=0}^{k-1} \varphi_j$ on $\mathbf{R}_{0+}$ with $\varphi_j(0) = x_j(0) = x^{(j)}(0) = x_{j0}$; $j \in \overline{k-1} \cup \{0\}$; $\forall t \in [-h, 0)$ for each given set of initial functions and $\phi(0) = \sum_{j=0}^{k-1} \frac{I_n}{j!}$ with $\varphi_j : [-h, 0] \to \mathbf{R}^n$, $j \in \overline{k-1} \cup \{0\}$ being bounded piecewise continuous with eventual discontinuities in a set of zero measure of $[-h, 0] \subset \mathbf{R}$ of bounded discontinuities, i.e. $\varphi_j \in BPC^{(0)}([-h, 0]), \mathbf{R}^n)$; $j \in \overline{k-1} \cup \{0\}$ and each given bounded piecewise continuous control $u : \mathbf{R}_{0+} \to \mathbf{R}^m$, with $u(t) = 0$ for $t \in [-h, 0)$, being a bounded piecewise continuous control function, and

(b)
$$x_\alpha(t) = \sum_{j=0}^{k-1} \left( \Phi_{\alpha j}(t) x_{j0} + \sum_{i=1}^{r} \int_{0}^{r_i} \Phi_\alpha(t-\tau) A_i \varphi_j(\tau - r_i) d\tau + \sum_{i=1}^{r} \int_{0}^{r_i} \Phi_\alpha(t-\tau) \tilde{A}_i(\tau) \varphi_j(\tau - r_i) d\tau \right)$$
$$+ \sum_{i=1}^{r} \int_{r_i}^{t} \Phi_\alpha(t-\tau) A_i x_\alpha(\tau - r_i) d\tau + \sum_{i=0}^{r} \int_{r_i}^{t} \Phi_\alpha(t-\tau) \tilde{A}_i(\tau) x_\alpha(\tau - r_i) d\tau$$
$$+ \int_0^t \Phi_\alpha(t-\tau) B(\tau) u(\tau) d\tau \; ; \; t \in \mathbf{R}_{0+} \tag{2.2}$$

which is time-differentiable satisfying (2.1) in $R_+$ with $k = [Re\,\alpha]+1$ if $\alpha\notin Z_+$ and $k = \alpha$ if $\alpha\in Z_+$, and

$$\Phi_{\alpha j}(t) := t^j E_{\alpha,j+1}(A_0 t^\alpha); \quad \Phi_\alpha(t) := t^{\alpha-1} E_{\alpha\alpha}(A_0 t^\alpha) \tag{2.3}$$

$$E_{\alpha j}(A_0 t^\alpha) = \sum_{\ell=0}^\infty \frac{(A_0 t^\alpha)^\ell}{\Gamma(\alpha\ell+j)}; \quad j\in\overline{k-1}\cup\{0,\alpha\} \tag{2.4}$$

for $t\in R_{0+}$ and $\Phi_{\alpha 0}(t) = \Phi_\alpha(t) = 0$ for $t<0$, where $E_{\alpha,j}(A_0 t^\alpha)$ are the Mittag-Leffler functions.

□

A technical result about norm upper-bounding functions of the matrix functions (2.3)-(2.4) follows:

**Lemma 2.2.** The following properties hold:

**(i)** There exist finite real constants $K_{E\alpha j}\geq 1$, $K_{\Phi\alpha j}\geq 1$; $j\in\overline{k-1}\cup\{0\}$ and $K_{\Phi\alpha}\geq 1$ such that for any $\alpha(\in R_+)<1$:

$$\left\|E_{\alpha j}(A_0 t^\alpha)\right\|\leq K_{E\alpha j}\left\|e^{A_0 t}\right\|, \left\|\Phi_{\alpha j}(t)\right\|\leq K_{\Phi\alpha j}\left\|t^j e^{A_0 t}\right\|; \; j\in\overline{k-1}\cup\{0,\alpha\}, \left\|\Phi_\alpha(t)\right\|\leq K_{\Phi\alpha}\left\|\frac{1}{t^{1-\alpha}}e^{A_0 t}\right\|$$

for $t(\in R_+)\geq 1$ \hfill (2.5)

**(ii)** If $\alpha(\in R_+)\geq 1$ then

$$\left\|E_{\alpha j}(A_0 t^\alpha)\right\| = \sum_{\ell=0}^\infty \left(\frac{\ell!}{\Gamma(\alpha\ell+j)}\right)\left(\frac{A_0^\ell t^{\alpha\ell}}{\ell!}\right) \leq \sup_{\ell\in Z_{0+}}\left(\frac{\ell!}{\Gamma(\ell+1)}\right)\left\|e^{A_0(t^\alpha)}\right\| = \left\|e^{A_0(t^\alpha)}\right\|; \; j\in\overline{k-1}\cup\{0\}$$

$$; \; t\in R_{0+}$$

$$\left\|\Phi_{\alpha j}(t)\right\| \leq \sup_{\ell\in Z_{0+}}\left(\frac{\ell!}{\Gamma(\alpha\ell+j+1)}\right) t^j\left\|e^{A_0(t^\alpha)}\right\| \leq t^j\left\|e^{A_0(t^\alpha)}\right\|, \; j\in\overline{k-1}\cup\{0\}; \; t\in R_{0+}$$

$$\left\|\Phi_\alpha(t)\right\| \leq \sup_{\ell\in Z_{0+}}\left(\frac{\ell!}{\Gamma((\ell+1)\alpha)}\right) t^{\alpha-1}\left\|e^{A_0(t^\alpha)}\right\| \leq t^{\alpha-1}\left\|e^{A_0(t^\alpha)}\right\|; \; t\in R_{0+} \tag{2.6}$$

If, in addition, $A_0$ is a stability matrix then $\left\|e^{A_0 t}\right\|\leq K e^{-\lambda t}$ and $\left\|e^{A_0(t^\alpha)}\right\|\leq K e^{-\lambda t^\alpha}\leq K e^{-\lambda t}$; $t\in R_{0+}$ for some real constants $K\geq 1$, $\lambda\in R_+$. Then, from (2.6)

$$\left\|E_{\alpha j}(A_0 t^\alpha)\right\|\leq K e^{-\lambda t}, \left\|\Phi_{\alpha j}(t)\right\|\leq t^j e^{-\lambda t}; \; j\in\overline{k-1}\cup\{0\}, \left\|\Phi_\alpha(t)\right\|\leq t^{\alpha-1} e^{-\lambda t} \tag{2.7}$$

for $t\in R_{0+}$ and the fractional dynamic system in the absence of delayed dynamics is exponentially stable if the standard fractional system for $\alpha = 1$ is exponentially stable.

**(iii)** The following inequalities hold:

$\left\|\Phi_{\alpha,k-1}(t)\right\|\leq t^{k-\alpha}\left\|\Phi_\alpha(t)\right\|$ for $\alpha\in(k-1,k]\cap R_+$ for $k\in Z_+$, $t\in R_{0+}$

$\left\|\Phi_\alpha(t)\right\|\leq t^{\alpha+1-k}\left\|\Phi_{\alpha,k-2}(t)\right\|$ for $\alpha\in[k-1,k)\cap R_+$, $t\in R_{0+}$

$\left\|\Phi_k(t)\right\|\equiv\left\|\Phi_{k,k-1}(t)\right\|$ for $\alpha = k\in Z_+$, $t\in R_{0+}$

**Proof**: Note from (2.3)-(2.4) for $0 < \alpha(\in \mathbf{R}_+) < 1$:

$$E_{\alpha j}(A_0 t^\alpha) = \sum_{\ell=0}^{\infty} \frac{(A_0 t^\alpha)^\ell}{\Gamma(\alpha \ell + j)} = \sum_{\ell=0}^{\infty} \left(\frac{\ell!}{\Gamma(\alpha \ell + j)}\right)\left(\frac{A_0^\ell t^{\alpha \ell}}{\ell!}\right) = \sum_{\ell=0}^{\infty} \left(\frac{t^{(\alpha-1)\ell}\ell!}{\Gamma(\alpha \ell + j)}\right)\left(\frac{A_0^\ell t^\ell}{\ell!}\right) \quad (2.8)$$

$$\|E_{\alpha j}(A_0 t^\alpha)\| \le \sup_{\tau \in (1,\infty) \cap \mathbf{R}} \left(\sup_{\ell \in \mathbf{Z}_{0+}} \left(\frac{\ell!}{\tau^{(1-\alpha)\ell}\Gamma(\alpha\ell+j)}\right)\right)\|e^{A_0 t}\|, j \in \overline{k-1} \cup \{0,\alpha\}; t \in (1,\infty) \cap \mathbf{R}_+ \quad (2.9)$$

$$\limsup_{t \to \infty} \|E_{\alpha,j}(A_0 t^\alpha)\| \le \limsup_{t \to \infty} \left(\sup_{\ell \in \mathbf{Z}_{0+}} \left(\frac{\ell!}{t^{(1-\alpha)\ell}\Gamma(\alpha\ell+j)}\right)\right)\left\|\sum_{\ell=0}^{\infty} \frac{A_0^\ell t^\ell}{\ell!}\right\| \le \limsup_{t \to \infty} \|e^{A_0 t}\|, j \in \overline{k-1} \cup \{0,\alpha\}$$

$$(2.10)$$

$$\|\Phi_{\alpha j}(t)\| \le \sup_{\tau \in (1,\infty) \cap \mathbf{R}} \left(\sup_{\ell \in \mathbf{Z}_{0+}} \left(\frac{\ell!}{\tau^{(1-\alpha)\ell}\Gamma(\alpha\ell+j+1)}\right)\right)\|t^j e^{A_0 t}\|, j \in \overline{k-1} \cup \{0,\alpha\}; t \in (1,\infty) \cap \mathbf{R}_+ \quad (2.11)$$

$$\limsup_{t \to \infty} \|\Phi_{\alpha j}(t)\| \le \limsup_{t \to \infty} \left(\sup_{\ell \in \mathbf{Z}_{0+}} \left(\frac{\ell!}{t^{(1-\alpha)\ell}\Gamma(\alpha\ell+j+1)}\right)\right)\left\|\sum_{\ell=0}^{\infty} \frac{A_0^\ell t^{\ell+j}}{\ell!}\right\| \le \limsup_{t \to \infty} \|t^j e^{A_0 t}\|, j \in \overline{k-1} \cup \{0\}$$

$$(2.12)$$

$$\|\Phi_\alpha(t)\| \le \sup_{\tau \in (1,\infty) \cap \mathbf{R}} \left(\sup_{\ell \in \mathbf{Z}_{0+}} \left(\frac{\ell!}{\tau^{(1-\alpha)\ell}\Gamma((\ell+1)\alpha)}\right)\right)\left\|\frac{1}{t^{1-\alpha}}e^{A_0 t}\right\|; t \in (1,\infty) \cap \mathbf{R}_+ \quad (2.13)$$

$$\limsup_{t \to \infty} \|\Phi_\alpha(t)\| \le \limsup_{t \to \infty} \left(\sup_{\ell \in \mathbf{Z}_{0+}} \left(\frac{\ell!}{t^{(1-\alpha)\ell}\Gamma((\ell+1)\alpha)}\right)\right)\left\|\sum_{\ell=0}^{\infty} \frac{A_0^\ell t^{\ell+\alpha-1}}{\ell!}\right\| \le \limsup_{t \to \infty} \left\|\frac{1}{t^{1-\alpha}}e^{A_0 t}\right\| \quad (2.14)$$

since

$$\limsup_{[\ell,\ell+1) \ni t(\in \mathbf{R}_{0+}) \to \infty} \left(\sup_{\ell \in \mathbf{Z}_{0+}} \left(\frac{\ell!}{t^{(1-\alpha)\ell}\Gamma((\ell+1)\alpha)}\right)\right) \le \limsup_{\mathbf{Z}_{0+} \ni \ell \to \infty} \left(\sup_{\ell \in \mathbf{Z}_{0+}} \left(\frac{\ell!}{\ell^{(1-\alpha)\ell}\Gamma(\alpha\ell+j+1)}\right)\right) = 0$$

$$(2.15)$$

The inequalities (2.5) hold since the above matrix norms are bounded on the real interval $(1,\infty)$ and their limit superior are upper-bounded by the given formulas and Property (i) is proved. On the other hand, if $\mathbf{R}_+ \ni \alpha \ge 1$ then

$$\|E_{\alpha j}(A_0 t^\alpha)\| = \sum_{\ell=0}^{\infty} \left(\frac{\ell!}{\Gamma(\alpha\ell+j)}\right)\left(\frac{A_0^\ell t^{\alpha\ell}}{\ell!}\right) \le \sup_{\ell \in \mathbf{Z}_{0+}}\left(\frac{\ell!}{\Gamma(\ell+1)}\right)\|e^{A_0(t^\alpha)}\| = \|e^{A_0(t^\alpha)}\|; j \in \overline{k-1} \cup \{0\}$$

$$; t \in \mathbf{R}_{0+} \quad (2.16)$$

$$\|\Phi_{\alpha j}(t)\| \le \sup_{\ell \in \mathbf{Z}_{0+}}\left(\frac{\ell!}{\Gamma(\alpha\ell+j+1)}\right)t^j \|e^{A_0(t^\alpha)}\| \le t^j \|e^{A_0(t^\alpha)}\|, j \in \overline{k-1} \cup \{0\}; t \in \mathbf{R}_{0+} \quad (2.17)$$

$$\|\Phi_\alpha(t)\| \le \sup_{\ell \in \mathbf{Z}_{0+}}\left(\frac{\ell!}{\Gamma((\ell+1)\alpha)}\right)t^{\alpha-1}\|e^{A_0(t^\alpha)}\| \le t^{\alpha-1}\|e^{A_0(t^\alpha)}\|; t \in \mathbf{R}_{0+} \quad (2.18)$$

If, in addition, $A_0$ is a stability matrix then $\|e^{A_0(t^\alpha)}\| \le Ke^{-\lambda t}$ and $\|e^{A_0(t^\alpha)}\| \le Ke^{-\lambda t^\alpha} \le Ke^{-\lambda t}$; $t \in \mathbf{R}_{0+}$ for some real constants $K \ge 1$ and $\lambda \in \mathbf{R}_+$ since $t^\alpha \ge t$; $\forall t \in \mathbf{R}_{0+}$. Properties (i)-(ii) have been proved.

**(iii)** It is proved as follows. Note from (2.3)-(2.4) that:

$$\|\Phi_{\alpha j}(t)\| = \left\|\sum_{\ell=0}^{\infty} \frac{A_0^\ell t^{\alpha\ell+j}}{\Gamma(\alpha\ell+j+1)}\right\| \leq \sup_{\ell \in Z_{0+}} \left(\frac{\Gamma((\ell+1)\alpha) t^{1+j-\alpha}}{\Gamma(\alpha\ell+j+1)}\right) \left\|\sum_{\ell=0}^{\infty} \frac{A_0^\ell t^{(\ell+1)\alpha-1}}{\Gamma((\ell+1)\alpha)}\right\|, t \in R_{0+}$$

so that if $k-1 < \alpha(\in R_+) \leq j+1 \leq k$. Then,

$$\|\Phi_{\alpha,k-1}(t)\| \leq \sup_{\ell \in Z_{0+}} \left(\frac{\Gamma((\ell+1)\alpha) t^{1+j-\alpha}}{\Gamma(\alpha\ell+k)}\right) \left\|\sum_{\ell=0}^{\infty} \frac{A_0^\ell t^{(\ell+1)\alpha-1}}{\Gamma((\ell+1)\alpha)}\right\| \leq t^{k-\alpha} \|\Phi_\alpha(t)\|, t \in R_{0+}$$

Also,

$$\|\Phi_\alpha(t)\| = \left\|\sum_{\ell=0}^{\infty} \frac{A_0^\ell t^{(\ell+1)\alpha-1}}{\Gamma((\ell+1)\alpha)}\right\| \leq t^{\alpha-j-1} \sup_{\ell \in Z_{0+}} \left(\frac{\Gamma(\alpha\ell+j+1)}{\Gamma((\ell+1)\alpha)}\right) \left\|\sum_{\ell=0}^{\infty} \frac{A_0^\ell t^{\alpha\ell+j}}{\Gamma(\alpha\ell+j+1)}\right\| \leq t^{\alpha-j-1} \|\Phi_{\alpha j}(t)\|, t \in R_{0+}$$

if $\sup_{\ell \in Z_{0+}} \left(\frac{\Gamma(\alpha\ell+j+1)}{\Gamma((\ell+1)\alpha)}\right) < \infty$. This implies that

$$\|\Phi_\alpha(t)\| \leq t^{\alpha+1-k} \|\Phi_{\alpha,k-2}(t)\| \text{ for } \alpha \in [k-1,k) \cap R_+$$

$$\|\Phi_k(t)\| \equiv \|\Phi_{k,k-1}(t)\| \text{ for } \alpha \in [k-1,k) \cap R_+ \qquad \square$$

## 3. Fixed point results

A technical definition is now given to facilitate the subsequent result about fixed point. Property (ii) has been proven.

**Definition 3.1.** $S(\overline{\varphi}, u)$ is the set of all the piecewise continuous n-vector function from $[-h, 0) \cup R_{0+}$ to $R^n$ being time-differentiable in $R_+$ which are solutions of (2.1) for all admissible k-tuples of initial conditions $\overline{\varphi} := (\varphi_0, \varphi_1, ..., \varphi_{k-1})$ with $\varphi_j \in BPC^{(0)}([-h, 0], R^n)$ and controls $u \in BPC^{(0)}(R_{0+}, R^n)$ with $\varphi_j(0) = x_j(0) = x^{(j)}(0) = x_{j0}; \forall j \in \overline{k-1} \cup \{0\}$. $\square$

A fixed point theorem is now given for the Caputo fractional system (2.1):

**Theorem 3.2.** Assume any set of r given finite delays $0 = r_0 < r_1 \leq .... \leq r_r = h < \infty$. The following properties hold:

**(i)** Assume that $\Phi_{\alpha j} \in L_\infty(R_{0+}, R^{n \times n})$ and $\int_0^\delta \|\Phi_\alpha(\delta-\tau)d\tau\| \|\widetilde{A}_0\|_\infty < 1$; and let $g_h: R_+ \to R_{0+}$ be defined by

$$g_h(\delta) := \left(1 - \int_0^\delta \|\Phi_\alpha(\delta-\tau)d\tau\| \|\widetilde{A}_0\|_\infty\right)^{-1} \times \left(\left\|\sum_{j=0}^{k-1} \Phi_{\alpha j}(\delta)\right\|\right.$$

$$\left. + \left(\left\|\int_0^\delta \Phi_\alpha(\delta-\tau)d\tau\right\|\right)\left(\sum_{i=1}^r \|\hat{A}_i\|_\infty\right)\right) \leq 1; \delta \in R_+ \qquad (3.1)$$

Then, the mapping $f_h : [-h, 0] \times R^n \to R_+ \times R^n$ defined by the state trajectory solution (2.2) of the uncontrolled system from any initial conditions in the admissible set is non-expansive and the solution is bounded fulfilling $\sup_{t \in R_{0+}} \|x_\alpha(t)\|_\infty \leq \sup_{t \in [-h,0]} \sum_{j=0}^{k-1} \left( \|\varphi_j(t)\|_\infty \right)$. If $g_h(\delta) \leq K_c(\delta) < 1$; $\forall \delta \in R_+$ then $f_h : [-h, 0] \times R^n \to R_+ \times R^n$ is contractive and possesses a unique fixed point, irrespective of the delays, in some bounded subset of $R^n$. Such a fixed point is $0 \in R^n$ which is a globally asymptotically stable equilibrium point.

**(ii)** Assume that $\Phi_{\alpha j} \in L_\infty(R_{0+}, R^{n \times n})$, $\Phi_\alpha \in L^2(R_{0+}, R^{n \times n})$ and $\int_0^\delta \|\Phi_\alpha(\delta - \tau) d\tau\| \left( \int_0^\delta \|\tilde{A}_0(t + \tau)\|^2 d\tau \right)^{1/2} < 1$; $\forall t \in R_{0+}$ and define $\hat{g}_h : R_{0+} \times R_+ \to R_{0+}$ point-wise as follows:

$$\hat{g}_h(t, \delta) := \left( 1 - \int_0^\delta \|\Phi_\alpha(\delta - \tau) d\tau\| \left( \int_0^\delta \|\tilde{A}_0(t + \tau)\|^2 d\tau \right)^{1/2} \right)^{-1} \times \left( \sum_{j=0}^{k-1} \|\Phi_{\alpha j}(\delta)\|_\infty + \left( \int_0^\delta \|\Phi_\alpha(\delta - \tau)\|^2 d\tau \right)^{1/2} \right.$$

$$\left. \times \left( \sum_{i=1}^r \left( \int_0^\delta \|\hat{A}_i(t + \tau - r_i)\|^2 \right)^{1/2} d\tau \right) \right) ; \quad \delta \in R_+ \tag{3.2}$$

Then, Property (i) still hold by replacing their corresponding constraints on $g_h$ by corresponding ones on $\hat{g}_h$.

**(iii)** Assume that a control $u(t) = \sum_{i=0}^r K_i(x_t, t) x(t - r_i)$ is injected to (2.1) where $K_i : R^n \times R_{0+} \to R^m$ is in $BPC(R_{0+}, R^m)$, $x_{it} : [\max(0, t - r_i), t] \to R^n$, $\forall i \in \overline{r-1} \cup \{0\}$; $\forall t \in R_{0+}$ is a strip of the state-trajectory solution of (2.1). Assume also that:

$$\|K_i(x_{it}, t)\|_\infty \leq K_i^0 < \infty, \forall i \in \overline{r-1} \cup \{0\}; \forall t \in R_{0+}, \Phi_{\alpha j} \in L_\infty(R_{0+}, R^{n \times n}), \Phi_\alpha \in L^1(R_{0+}, R^{n \times n})$$

(3.3)

and define $\hat{g}_f : R_+ \to R_{0+}$ as

$$g_h(\delta) := \left( 1 - \int_0^\delta \|\Phi_\alpha(\delta - \tau) d\tau\| \left( \|\tilde{A}_0\|_\infty + \|B\|_\infty K_0^0 \right) \right)^{-1} \times \left( \left\| \sum_{j=0}^{k-1} \Phi_{\alpha j}(\delta) \right\| \right.$$

$$\left. + \left( \left\| \int_0^\delta \Phi_\alpha(\delta - \tau) d\tau \right\| \right) \left( \sum_{i=1}^r \left( \|\hat{A}_i\|_\infty + \|B\|_\infty K_i^0 \right) \right) \right) \leq 1; \quad \delta \in R_+ \tag{3.4}$$

provided that $\int_0^\delta \|\Phi_\alpha(\delta - \tau) d\tau\| \left( \|\tilde{A}_0\|_\infty + \|B\|_\infty K_0^0 \right) < 1$. Then, for any given set of finite delays, the mapping $f_f : [-h, 0] \times R^n \times R^m \times R_{0+} \to R_+ \times R^n$ defined by the state trajectory solution (2.2) of the controlled system from any initial conditions in the admissible set and any given admissible control is a non-expansive mapping if $g_h(\delta) \leq 1$; $\forall \delta \in R_+$ and contractive and the zero equilibrium is the unique

fixed point, irrespective of the delays and control, if $g_h(\delta) \leq K_c(\delta) < 1; \forall \delta \in \mathbf{R}_+$ which is also a globally asymptotically stable equilibrium point.

**(iv)** Assume that $\exists \varepsilon (<1) \in \mathbf{R}_+$ that $g_h(\delta) < 1 - \varepsilon$; $\forall t \in \mathbf{R}_{0+}$. Then, state trajectory solution (2.2) of the forced system from any initial conditions in the admissible set is defined by contractive self-mapping with a unique fixed point in some bounded subset of $\mathbf{R}^n$ for all controls of the form $u(t) = \sum_{i=0}^{r} K_i(x_{it}, t) x(t - r_i)$ fulfilling $\| K_i(x_{it}, t) \|_\infty \leq \dfrac{\varepsilon}{(r+1)\left( \left\| \int_0^\delta \Phi_\alpha(\delta - \tau) d\tau \right\| \right) \| B \|_\infty}$, $\forall i \in \overline{r-1} \cup \{0\}$

**(v)** Assume that $\Phi_{\alpha j} \in L_\infty(\mathbf{R}_{0+}, \mathbf{R}^{n \times n})$; $\forall j \in \overline{k-1} \cup \{0\}$, $\Phi_\alpha \in L^2(\mathbf{R}_{0+}, \mathbf{R}^{n \times n})$ and $B K_i \in L^2(\mathbf{R}_{0+}, \mathbf{R}^{n \times n})$; $\forall i \in \overline{r-1} \cup \{0\}$, instead of the hypotheses (3.3), and define $\hat{g}_f : \mathbf{R}_{0+} \times \mathbf{R}_+ \to \mathbf{R}_{0+}$ as:

$$\hat{g}_f(t,\delta) := \left(1 - \int_0^\delta \| \Phi_\alpha(\delta - \tau) d\tau \| \left( \left( \int_0^\delta \| \tilde{A}_0(t+\tau) \|^2 d\tau \right)^{1/2} + \left( \int_0^\delta \| B(t+\tau) K_0(x_{t+\tau}, t+\tau) \|^2 d\tau \right)^{1/2} \right) \right)^{-1} \times \left( \sum_{j=0}^{k-1} \| \Phi_{\alpha j}(\delta) \|_\infty \right.$$
$$\left. + \left( \int_0^\delta \| \Phi_\alpha(\delta - \tau) \|^2 d\tau \right)^{1/2} \left( \sum_{i=1}^{r} \left( \int_0^\delta \| \hat{A}_i(t+\tau - r_i) \|^2 \right)^{1/2} d\tau + \left( \int_0^\delta \| B(t+\tau) K_i(x_{t+\tau}, t+\tau) \|^2 d\tau \right)^{1/2} \right) \right)$$

$$; \forall t, \delta \in \mathbf{R}_+ \quad (3.5)$$

Provided that the inverse exists on $\mathbf{R}_{0+}$, Then, Property (iii) still hold by replacing their corresponding constraints on $g_f$ by corresponding ones on $\hat{g}_f$. If, in addition, $\exists \varepsilon(<1) \in \mathbf{R}_+$, $\delta = \delta(\varepsilon) \in \mathbf{R}_+$ such that $\hat{g}_h(\delta) < 1 - \varepsilon$; $\forall t \in \mathbf{R}_{0+}$ then the mapping $f_f : [-h, 0] \times \mathbf{R}^n \times \mathbf{R}^m \times \mathbf{R}_{0+} \to \mathbf{R}_+ \times \mathbf{R}^n$ defining the state-trajectory solution from any set of admissible initial conditions and all controls $u(t) = \sum_{i=0}^{r} K_i(x_{it}, t) x(t - r_i)$ being subject to

$$\sum_{i=1}^{r} \int_0^\delta \left( \| B(t+\tau) K_i(x_{t+\tau}, t+\tau) \|^2 d\tau \right)^{1/2} \leq \dfrac{\varepsilon}{\left( \int_0^\delta \| \Phi_\alpha(\delta - \tau) \|^2 d\tau \right)^{1/2} \| B \|_\infty}; \forall t \in \mathbf{R}_+ \quad (3.6)$$

is contractive with a unique fixed point, irrespective of the delays, which is $0 \in \mathbf{R}^n$ bein a globally asymptotically stable equilibrium point

**Proof.** The point-wise difference between two solutions x(t) and z(t) of (2.1) subject to respective piecewise continuous initial conditions $\varphi_x : [-h, 0] \to \mathbf{R}^n$ and $\varphi_z : [-h, 0] \to \mathbf{R}^n$ and respective controls $u_x, u_y \in BPC^{(0)}(\mathbf{R}_{0+}, \mathbf{R}^n)$ is according to (2.2):

$$x_\alpha(t) - z_\alpha(t) = \sum_{j=0}^{k-1} \left( \Phi_{\alpha j}(t)(x_{j0} - z_{j0}) + \sum_{i=1}^{r} \int_0^{r_i} \Phi_\alpha(t-\tau) A_i (\varphi_{xj}(\tau - r_i) - \varphi_{zj}(\tau - r_i)) d\tau \right)$$
$$+ \sum_{j=0}^{k-1} \sum_{i=1}^{r} \int_0^{r_i} \Phi_\alpha(t-\tau) \tilde{A}_i(\tau) (\varphi_{xj}(\tau - r_i) - \varphi_{zj}(\tau - r_i)) d\tau$$

$$+\sum_{i=1}^{r}\int_{r_i}^{t}\Phi_\alpha(t-\tau)A_i(x_\alpha(\tau-r_i)-z_\alpha(\tau-r_i))d\tau+\sum_{i=0}^{r}\int_{0}^{r_i}\Phi_\alpha(t-\tau)\tilde{A}_i(\tau)(x_\alpha(\tau-r_i)-z_\alpha(\tau-r_i))d\tau$$

$$+\int_{0}^{t}\Phi_\alpha(t-\tau)B(\tau)(u_x(\tau)-u_z(\tau))d\tau\; ;\; t\in R_{0+} \quad (3.7)$$

Note from (2.3) that $\Phi_{\alpha j}(0)=\frac{I_n}{j!}$; $\forall j\in \overline{k-1}\cup\{0\}$ what is used in the definition of the metric space $(M,\|.\|_\infty)$ with the supremum metric $\|.\|_\infty$

$$M:=\left\{\phi\in PBC^{(0)}([-h,0)\cup R_{0+},R^n): \phi\in S(\bar{\phi},u), \phi\equiv\left(\sum_{j=0}^{k-1}\phi_j\right)\in BPC^{(0)}([-h,0),R^n)\right.$$

$$\left. ;\forall j\in\overline{k-1}\cup\{0\}, \phi_u\in M_u\right\} \quad (3.8)$$

where

$$M_u:=\left\{\phi\in PBC^{(0)}(R_{0+},R^n):\phi\in S(0,u), u\in BPC^{(0)}(R_{0+},R^n)\right\} \quad (3.9)$$

where $PBC^{(0)}(R,R^n)$ is the set of bounded continuous n-vector functions on $R$. Now, define $P:M\to M$ as the subsequent piecewise bounded continuous function on $[-h,0)\cup R_{0+}$, which is bounded continuous on $R_+$, i.e. $\phi\in PBC^{(0)}([-h,0]R^n)$, $\phi\in BC^{(0)}(R_+,R^n)$ and satisfies (2.1) on $R_+$. One gets for any bounded piecewise continuous solution of (2.1):

$$(P(\phi,u))(t):=\sum_{j=0}^{k-1}\Phi_{\alpha j}(t)\phi_j(0)+\int_0^t\Phi_\alpha(t-\tau)\tilde{A}_0(\tau)\phi(\tau)d\tau$$

$$+\sum_{j=0}^{k-1}\sum_{i=1}^{r}\left(\int_0^{r_i}\Phi_\alpha(t-\tau)\hat{A}_i(\tau)\phi_j(\tau-r_i)d\tau+\int_{r_i}^{t}\Phi_\alpha(t-\tau)\hat{A}_i(\tau)\phi_j(\tau-r_i)d\tau\right)$$

$$+\int_0^t\Phi_\alpha(t-\tau)B(\tau)u_\phi(\tau)d\tau \quad (3.10)$$

Note that the supremum metric $[-h,0)\cup R_{0+}$ is induced by the supremum norm on $[-h,0)\cup R_{0+}$ so that it is then the supremum norm. Define the truncated $\phi_{a,b}\in M$ as $\phi_{a,b}(\tau)=\phi(\tau)$, $\tau\in[a,b)$ and $\phi_{a,b}(\tau)=0$, $\tau\in[0,a)\cup[b,\infty)\subset R_{0+}; \forall a,b(>a)\in[a,b)\subseteq R_{0+}$, $\forall t\in R_{0+}$, $\forall \phi\in M$ and note that $\|\phi\|_{a,b}=\|\phi_{a,b}\|_\infty\le\|\phi\|_\infty$ where $\|\phi\|_{a,b}:=\sup_{\tau\in[a,b)}\|\phi(\tau)\|$. $\phi_t\in M$ is the truncated $\phi\in M$ on $[t-h,t)$. Norms without subscripts mean vector or correspondingly induced matrix norms (as, for instance, the $\ell_2$- vector or induced matrix norms) or point-wise values of such norms for vector or matrix functions in the subsequent developments. Let $M_t$ be the space of truncated functions $\phi_t\in M$. Note that any truncated solution of (2.1) on any finite interval is always in $M$ so that one gets from (3.8) for any $\delta\in R_+$ in the most general controlled case with control $u(t)=\sum_{i=0}^{r}K_i(x_t,t)x(t-r_i)$:

$$\left\|\left(P(\phi,u_\phi)\right)(t+\delta)-\left(P(\eta,u_\eta)\right)(t+\delta)\right\| \leq \sum_{j=0}^{k-1}\left\|\Phi_{\alpha j}(\delta)\right\|\left\|\phi_j(t)-\eta_j(t)\right\|_t$$

$$+\left\|\int_0^\delta \Phi_\alpha(\delta-\tau)B(t+\tau)d\tau\right\|\left\|u_\phi-u_\eta\right\|_{t+\delta}$$

$$+\left(\sum_{i=1}^{r}\left\|\int_0^\delta \Phi_\alpha(\delta-\tau)\hat{A}_i(t+\tau)d\tau\right\|_{t+\delta-r_i}+\int_0^\delta\left\|\Phi_\alpha(\delta-\tau)\tilde{A}_0(t+\tau)d\tau\right\|\right)\|\phi-\eta\|_{t+\delta} \quad (3.11)$$

$$\leq \left(\left\|\sum_{j=0}^{k-1}\Phi_{\alpha j}(\delta)\right\|+\sum_{i=1}^{r}\int_0^\delta\left\|\Phi_\alpha(\delta-\tau)B(t+\tau)K_i\left(x_{i,t+\tau},t+\tau\right)d\tau\right\|\|\phi-\eta\|_{t+\delta-r_i}\right.$$

$$+\sum_{i=1}^{r}\int_0^\delta\left\|\Phi_\alpha(\delta-\tau)\hat{A}_i(t+\tau)d\tau\right\|\|\phi-\eta\|_{t+\delta-r_i}$$

$$\left.+\int_0^\delta\left\|\Phi_\alpha(\delta-\tau)\left(\tilde{A}_0(t+\tau)+B(t+\tau)K_0(t+\tau)\right)\right\|d\tau\|\phi-\eta\|_{t+\delta}\right) \quad (3.12)$$

$$\leq \left(\left\|\sum_{j=0}^{k-1}\Phi_{\alpha j}(\delta)\right\|+\sum_{i=1}^{r}\int_0^\delta\|\Phi_\alpha(\delta-\tau)d\tau\|\left(\sum_{i=1}^{r}\left(\|\hat{A}_i\|_\infty+\|B\|_\infty K_i^0\right)\|\phi-\eta\|_{t+\delta-r_i}+\left(\|\tilde{A}_0\|_\infty+\|B\|_\infty K_0^0\right)\|\phi-\eta\|_{t+\delta}\right)\right)$$

$$\leq \left(\left\|\sum_{j=0}^{k-1}\Phi_{\alpha j}(\delta)\right\|+\sum_{i=1}^{r}\int_0^\delta\|\Phi_\alpha(\delta-\tau)d\tau\|\left(\sum_{i=1}^{r}\left(\|\hat{A}_i\|_\infty+\|B\|_\infty K_i^0\right)+\left(\|\tilde{A}_0\|_\infty+\|B\|_\infty K_0^0\right)\right)\right)\|\phi-\eta\|_{t+\delta} \quad (3.13)$$

where the property that $A_0$ is constant has been used to rewrite the limits of the involved integral is the most convenient fashion to simplify the related expressions. Eq. (3.13a) leads to:

$$\|\phi-\eta\|_{t+\delta} \leq \left(1-\int_0^\delta\|\Phi_\alpha(\delta-\tau)d\tau\|\left(\|\tilde{A}_0\|_\infty+\|B\|_\infty\|K_0^0\|_\infty\right)\right)^{-1}$$

$$\times\left(\left\|\sum_{j=0}^{k-1}\Phi_{\alpha j}(\delta)\right\|+\sum_{i=1}^{r}\int_0^\delta\|\Phi_\alpha(\delta-\tau)d\tau\|\left(\sum_{i=1}^{r}\left(\|\hat{A}_i\|_\infty+\|B\|_\infty\|K_i^0\|_\infty\right)\|\phi-\eta\|_{t+\delta-r_i}\right)\right)$$

$$\leq \left(1-\int_0^\delta\|\Phi_\alpha(\delta-\tau)d\tau\|\left(\|\tilde{A}_0\|_\infty+\|B\|_\infty\|K_0^0\|_\infty\right)\right)^{-1}$$

$$\times\left(\left\|\sum_{j=0}^{k-1}\Phi_{\alpha j}(\delta)\right\|+\sum_{i=1}^{r}\int_0^\delta\|\Phi_\alpha(\delta-\tau)d\tau\|\left(\sum_{i=1}^{r}\left(\|\hat{A}_i\|_\infty+\|B\|_\infty\|K_i^0\|_\infty\right)\|\phi-\eta\|_{t+\delta-r_1}\right)\right)$$

$$;\forall\delta\in\mathbf{R}_+,\forall t\in\mathbf{R}_{0+} \quad (3.14)$$

provided that $\int_0^\delta\|\Phi_\alpha(\delta-\tau)d\tau\|\left(\|\tilde{A}_0\|_\infty+\|B\|_\infty K_0^0\right)<1$, since $r_1\leq r_i\,(i\in\bar{r})$, so that

$$\|\phi-\eta\|_t \leq \left(1-\int_0^\delta\|\Phi_\alpha(\delta-\tau)d\tau\|\left(\|\tilde{A}_0\|_\infty+\|B\|_\infty K_0^0\right)\right)^{-1}$$

$$\times\left(\left\|\sum_{j=0}^{k-1}\Phi_{\alpha j}(\delta)\right\|+\sum_{i=1}^{r}\int_0^\delta\|\Phi_\alpha(\delta-\tau)d\tau\|\left(\sum_{i=1}^{r}\left(\|\hat{A}_i\|_\infty+\|B\|_\infty\|K_i^0\|_\infty\right)\|\phi-\eta\|_{t-r_1}\right)\right)$$

$$= g_h(\delta)\|\phi-\eta\|_{t-r_1};\;\forall\delta\in\mathbf{R}_+,\forall t\in\mathbf{R}_{0+} \quad (3.15)$$

Then, the mapping $f_h : [-h, 0] \times R^n \to R_+ \times R^n$ defining the state trajectory solution from admissible initial conditions is non-expansive if $g_h(\delta) \leq 1$. Furthermore, the state trajectory solution is globally Lyapunov stable since by taking the trivial solution $\eta \equiv 0$ in (3.15), it follows that any solution $\phi$ of (2.1) generated from any set of admissible initial conditions is uniformly bounded on $R_{0+}$. If, in addition, $g_h(\delta) \leq K_c(\delta) < 1$ then it follows also from (3.15) that any real sequence of the form $\{v(kr_1 + \tau)\}_{k \in Z_{0+}, \tau \in [0, r_1) \cap R_{0+}}$ is a convergent Cauchy sequence so that the metric space the metric space $(M, \|.\|_\infty)$ of the solutions of (2.1) under the class of given initial conditions and controls with the supremum metric $\|.\|_\infty$ is complete. Therefore, a unique fixed point exists on some bounded set of $R^n$ from Banach contraction principle. Since

$$\lim_{Z_{0+} \ni k \to \infty, \tau \in [0, r_1) \cap R_{0+}} \|\phi - \eta\|_{(k+1)(r_1+\tau)} \leq \left(\lim_{Z_{0+} \ni k \to \infty} K_c^k(\delta)\right) \|\phi - \eta\|_{k(r_1+\tau), \tau \in [0, r_1) \cap R_{0+}} = 0$$

(3.16)

it follows by taking one of the solutions to be the trivial solution that the only fixed point is the equilibrium point zero which is a globally asymptotically stable attractor. Property (i) has been proved. By zeroing the control and considering the uncontrolled system, one proves Property (i) as a particular case of property (iii). Property (ii) (and its particular case Property (iv) for the case of controller gains satisfying $\|K_i(x_{it}, t)\|_\infty \leq K_i^0 < \infty$) and Property (v) are proved by using similar technical tools to those involved in the above proofs by replacing the basic inequality (3.13) by

$$\|(P(\phi, u_\phi))(t) - (P(\eta, u_\eta))(t)\| \leq \left(\left\|\sum_{j=0}^{k-1} \Phi_{\alpha j}(t)\right\|_{t-r_1, t} + \sum_{i=1}^{r}\left(\int_{t-r_i}^{t} \|\Phi_\alpha(t-\tau)\|^2 d\tau\right)^{1/2}\right)$$
$$\times \left(\left(\int_{t-r_i}^{t} \|\hat{A}_i(\tau)\|^2\right)^{1/2} d\tau + \int_{t-r_i}^{t} \left(\|B(\tau)K_i(x_\tau, \tau)\|^2 d\tau\right)^{1/2}\right) + \left(\int_{t-r_1}^{t} \|\tilde{A}_0(\tau)\|^2 d\tau\right)^{1/2} \|\phi - \eta\|_{t-r_1}$$
$$; \forall \delta \in R_+, \forall t \in R_{0+} \quad (3.17)$$

□

If all the delays are zero, it is more convenient to discuss the ad-hoc solution version of (2.2):

$$x_\alpha(t) = \sum_{j=0}^{k-1} \overline{\Phi}_{\alpha j}(t) x_{j0} + \sum_{i=0}^{r} \int_0^t \overline{\Phi}_\alpha(t-\tau) \tilde{A}_i(\tau) x_\alpha(\tau) d\tau + \int_0^t \overline{\Phi}_\alpha(t-\tau) B(\tau) u(\tau) d\tau \, ; \, t \in R_{0+}$$

(3.18)

where $\overline{\Phi}_{\alpha j}(t)$ and $\overline{\Phi}_\alpha(t)$ are similar to $\Phi_{\alpha j}(t)$ from Eqs.(2.3)-(2.4) by replacing $A_0 \to \left(\sum_{i=0}^{r} A_i\right)$

The following result is a parallel one to Theorem 3.2 for the case of absence of delays:

**Theorem 3.3**. Assume that:

(1) $\|\overline{\Phi}_{\alpha j}(t)\| \leq K_{0j}(t)$; $t \in \boldsymbol{R}_{0+}$, $\forall j \in \overline{k-1} \cup \{0\}$ with $\max_{0 \leq j \leq k-1}\left(\sup_{t \in \boldsymbol{R}_{0+}} K_{0j}(t)\right) \leq \overline{K}_0 < \infty$

(2) $\overline{\Phi}_\alpha \in L^1(\boldsymbol{R}_{0+}, \boldsymbol{R}^{n \times n})$ with $\sup_{t \in \boldsymbol{R}_{0+}} \|\overline{\Phi}_\alpha(t)\| \leq \overline{K}_1 < \infty$

Then, the Caputo delay-free fractional dynamic system (2.1) of real order $\alpha$ has the following properties:

**(i)** It is globally stable under a control $u(t) = \sum_{i=0}^{r} K_i(x_{it}, t) x(t)$ subject to $\|K_i(x_{it}, t)\|_\infty \leq K_i^0 < \infty$, $\forall i \in \overline{r-1} \cup \{0\}$ if $\overline{K}_1 < \dfrac{1}{\sum_{i=0}^{r}\left(\|\tilde{A}_i\|_\infty + \|B\|_\infty K_i^0\right)}$, $\forall i \in \overline{r-1} \cup \{0\}$. If, in addition, $K_{0j}(t) \to 0$ as $t \to \infty$; $\forall j \in \overline{k-1} \cup \{0\}$ then the system is globally asymptotically stable to the zero equilibrium point.

**(ii)** Property (i) holds if $K_0(t) = K_0^0$ is constant if $\overline{K}_1 < \dfrac{1}{\sum_{i=0}^{r}\|\tilde{A}_i\|_\infty + \|B\|_\infty \left(\sum_{i=1}^{r} K_i^0\right)}$ where $\overline{\Phi}_{\alpha j}(t)$ and $\overline{\Phi}_\alpha(t)$ are similar to $\Phi_{\alpha j}(t)$ from Eqs.(2.3)-(2.4) by replacing $A_0 \to \left(\sum_{i=0}^{r} A_i\right) + BK_0^0$.

**Proof**: **(i)** One gets after taking norms in (3.18) that:

$$\|x_\alpha(t)\| \leq \left\|\sum_{j=0}^{k-1} \overline{\Phi}_{\alpha j}(t) x_{j0}\right\| + \sum_{i=0}^{r}\left(\|\tilde{A}_i\|_\infty + \|B\|K_i^0\right)\left(\int_0^t \|\Phi_\alpha(t-\tau)\| d\tau\right) \sup_{\tau \in [0,t)} \|x_\alpha(\tau)\|$$

$$\leq \sum_{j=0}^{k-1} K_{0j}(t)\|x_{j0}\| + \sum_{i=0}^{r}\left(\|\tilde{A}_i\|_\infty + \|B\|K_i^0\right)\left(\int_0^t \|\Phi_\alpha(t-\tau)\| d\tau\right) \sup_{\tau \in [0,t)} \|x_\alpha(\tau)\|$$

$$\leq \sum_{j=0}^{k-1} K_{0j}(t)\|x_{j0}\| + \overline{K}_1 \left(\sum_{i=0}^{r}\left(\|\tilde{A}_i\|_\infty + \|B\|K_i^0\right)\right) \sup_{\tau \in [0,t)} \|x_\alpha(\tau)\| \qquad (3.19)$$

$$\leq \overline{K}_0 \left(\sum_{j=0}^{k-1}\|x_{j0}\|\right) + \overline{K}_1 \left(\sum_{i=0}^{r}\left(\|\tilde{A}_i\|_\infty + \|B\|K_i^0\right)\right) \sup_{\tau \in [0,t)} \|x_\alpha(\tau)\| \; ; \; t \in \boldsymbol{R}_{0+} \qquad (3.20)$$

Then, one gets from (3.20):

$$\|x_\alpha(t)\| \leq \sup_{\tau \in [0,t)}\|x_\alpha(\tau)\| \leq \left(1 - \overline{K}_1\left(\sum_{i=0}^{r}\left(\|\tilde{A}_i\|_\infty + \|B\|K_i^0\right)\right)\right)^{-1} \overline{K}_0 \left(\sum_{j=0}^{k-1}\|x_{j0}\|\right) \leq \overline{K}_2 < \infty$$

$$; t \in \boldsymbol{R}_{0+} \quad (3.21)$$

where $\overline{K}_2 := \sup_{t \in \boldsymbol{R}_{0+}} \|x_\alpha(t)\| < \infty$. As a result, the Caputo fractional system of real order $\alpha$ is globally stable under zero delays since any state trajectory solution generated from any admissible initial conditions is bounded for all time. The proof of property (ii) is similar to that of (i) under the modified constraints. Now, assume that if, in addition, $K_{0j}(t) \to 0$ as $t \to \infty$; $\forall j \in \overline{k-1} \cup \{0\}$, then

$$\|x_\alpha(t)\| \leq \min\left(1, \left[\sum_{j=0}^{k-1} K_{0j}(t) + \sum_{i=0}^{r}\left(\|\widetilde{A}_i\|_\infty + \|B\|K_i^0\right)\overline{K}_1\right]\right)\overline{K}_2 \qquad (3.21)$$

so that

$$\limsup_{t\to\infty}\|x_\alpha(t)\| \leq \min\left(1, \limsup_{t\to\infty}\left[\sum_{j=0}^{k-1} K_{0j}(t) + \sum_{i=0}^{r}\left(\|\widetilde{A}_i\|_\infty + \|B\|K_i^0\right)\overline{K}_1\right]\right)\overline{K}_2$$

$$= \left(\sum_{i=0}^{r}\left(\|\widetilde{A}_i\|_\infty + \|B\|K_i^0\right)\overline{K}_1\right)\overline{K}_2 < \overline{K}_2 \qquad (3.22)$$

since $\lim_{t\to\infty} K_{0j}(t) = 0$; $\forall j \in \overline{k-1} \cup \{0\}$ and $\sum_{i=0}^{r}\left(\|\widetilde{A}_i\|_\infty + \|B\|K_i^0\right)\overline{K}_1 < 1$. Eq. (3.22) implies that the supremum $\|x_\alpha(t)\|$ on $\mathbf{R}_{0+}$ is reached by the first time at some finite time $t_0 \in \mathbf{R}_{0+}$. Thus, one gets from (3.19) that

$$\lim_{t\to\infty}\|x_\alpha(t)\| \leq \lim_{t\to\infty}\left(\sup_{\tau\in[t_0,t]}\|x_\alpha(\tau)\|\right) \leq \left(1 - \overline{K}_1\left(\sum_{i=0}^{r}\left(\|\widetilde{A}_i\|_\infty + \|B\|K_i^0\right)\right)\right)^{-1}\left(\sum_{j=0}^{k-1}\|x_{\alpha j t_0}\|\left(\lim_{t\to\infty}\left(K_{0j}(t-t_0)\right)\right)\right) = 0$$

(3.23)

provided that $K_{0j}(t) \to 0$ as $t \to \infty$; $\forall j \in \overline{k-1} \cup \{0\}$ which proves the global asymptotic stability. Property (i) has been proved. Property (ii) follows in a similar way under the modified constraints $K_0(t) = K_0^0$ $\overline{K}_1 < \dfrac{1}{\sum_{i=0}^{r}\|\widetilde{A}_i\|_\infty + \|B\|_\infty\left(\sum_{i=1}^{r}K_i^0\right)}$, $\overline{\Phi}_{\alpha j}(t)$, $\overline{\Phi}_\alpha(t)$ being similar to $\Phi_{\alpha j}(t)$ from Eqs.(2.3)-(2.4) by replacing $A_0 \to \left(\sum_{i=0}^{r} A_i\right) + BK_0^0$. $\square$

The subsequent stability result is based on a transformation of the matrix $A_0$ to its diagonal Jordan form which allows an easy computation of the $\ell_2$ - matrix measure of its diagonal part:

**Theorem 3.4.** Assume that $J_{A_0} = J_{A_{0d}} + \widetilde{J}_{A_0}$ is the Jordan form of $A_0$ with $J_{A_{0d}}$ being diagonal and $\widetilde{J}_{A_0}$ being off-diagonal such that the above decomposition is unique with $A_0 = T^{-1}J_{A_0}T$ where T is a unique non-singular transformation matrix. The following properties hold:

(i) The Caputo fractional differential system (2.1) is globally Lyapunov stable independently if the $\ell_2$ - matrix measure of $J_{A_{0d}}$ is negative, i.e. $\mu_2\left(J_{A_{0d}}^{1/\alpha}\right) = \dfrac{1}{2}\lambda_{max}\left(J_{A_{0d}}^{1/\alpha} + J_{A_{0d}}^{*1/\alpha}\right) = \max_{k\in\overline{n}} Re\left(\lambda_k^{1/\alpha}\right) < 0$;

$\forall \lambda_k \in \sigma\left(J_{A_{0d}}^{1/\alpha}\right)$, that is the spectrum of $J_{A_{0d}}^{1/\alpha}$ and, furthermore,

$$\left\|\left(\dfrac{1}{\beta_0}T^{-1}\widetilde{J}_{A_0}T, \dfrac{1}{\beta_1}T^{-1}A_1T, \dfrac{1}{\beta_2}T^{-1}A_2T, \cdots, \dfrac{1}{\beta_r}T^{-1}A_rT\right)\right\|_2 \leq \left|\mu_2\left(J_{A_{0d}}\right)\right|^{1/\alpha} \qquad (3.24)$$

for some set of numbers $\beta_i \in \mathbf{R}_+$ $(i \in \bar{p} \cup \{0\})$ satisfying $\sum_{i=0}^{r} \beta_i^2 = 1$. The fractional system is globally asymptotically Lyapunov stable for one such a set of real numbers if $\mu_2\left(J_{A_{0d}}^{1/\alpha}\right) < 0$, and

$$\left\|\left(\frac{1}{\beta_0} T^{-1} \tilde{J}_{A_0} T, \frac{1}{\beta_1} T^{-1} A_1 T, \frac{1}{\beta_2} T^{-1} A_2 T, \cdots, \frac{1}{\beta_r} T^{-1} A_r T\right)\right\|_2 < \left|\mu_2\left(J_{A_{0d}}\right)\right|^{1/\alpha} \quad (3.25)$$

**(ii)** A necessary condition for $\mu_2\left(J_{A_{0d}}^{1/\alpha}\right) < 0$, being a necessary condition for testing (3.25), is that $A_0$ be a stability matrix with $|arg(\lambda)| < \frac{\alpha \pi}{2}$; $\forall \lambda \in \sigma(A_0)$. Such a condition holds directly if $\alpha > 2\varphi/\pi$ where $(-\varphi, \varphi) \subseteq \left(-\frac{\pi}{2}, \frac{\pi}{2}\right)$ is the symmetric maximum real interval containing the arguments of all $\lambda \in \sigma(A_0)$. It also holds, in particular, if $A_0$ is a stability matrix and $\alpha (\in \mathbf{R}_+) \geq 1$.

**Proof**: It follows by using the matrix similarity transformation $A_0 = T^{-1} J_{A_0} T = T^{-1}\left(J_{A_{0d}} + \tilde{J}_{A_0}\right) T$ and using the homogeneous transformed Caputo fractional differential system from (2.1):

$$\left({}^C D_{0+}^{\alpha} z\right)(t) = \left({}^C D_{0+}^{\alpha} T x\right)(t) = \sum_{i=0}^{r} A_i T x(t - h_i) \Leftrightarrow$$

$$\left({}^C D_{0+}^{\alpha} x\right)(t) = \sum_{i=0}^{p} T^{-1} A_i T x(t - h_i) = T^{-1} A_0 T x(t) + \sum_{i=1}^{r} T^{-1} A_i T x(t - h_i)$$

$$= T^{-1} J_{A_{0d}} T x(t) + \sum_{i=0}^{r} T^{-1} \overline{A}_i T x(t - h_i) \quad (3.26)$$

where $z(t) = T x(t)$; $\forall t \in \mathbf{R}_{0+}$, $h_0 = 0$ plays the role of an additional delay. $\overline{A}_0 = \tilde{J}_{A_0}$ and $\overline{A}_i = A_i$ ($i \in \bar{r}$) by noting also that since $\left(J_{A_{0d}} + J_{A_{0d}}^*\right)$ is diagonal with real eigenvalues by construction, one has:

$$\left|\mu_2\left(J_{A_{0d}}^{1/\alpha}\right)\right| = \left|\frac{1}{2} \lambda_{max}\left(J_{A_{0d}}^{1/\alpha} + J_{A_{0d}}^{1/\alpha *}\right)\right| = \left|\lambda_{max}\left(J_{A_{0d}}^{1/\alpha}\right)\right|$$

$$= \left|Re\, \lambda_{max}\left(J_{A_{0d}}^{1/\alpha}\right)\right| = \left|Re\, \lambda_{max}^{1/\alpha}\left(J_{A_{0d}}\right)\right| = \left|Re\, \lambda_{max}^{1/\alpha}(A_{0d})\right| = \left|\mu_2\left(J_{A_{0d}}\right)\right|^{1/\alpha} \quad (3.27)$$

Then, the remaining part of the proof of property (i) is similar to that quoted as a sufficient condition for stability independent of the delays in [27]. Property (i) has been proved.

To prove property (ii), note that:

$$\lambda_{max}\left(J_{A_{0d}}^{1/\alpha} + J_{A_{0d}}^{*\,1/\alpha}\right) \neq \frac{1}{2} \lambda_{max}\left(J_{A_{0d}} + J_{A_{0d}}^*\right)^{1/\alpha} \text{ if } \mu_2\left(J_{A_{0d}}^{1/\alpha}\right) < 0 \text{ so that}$$

$$0 > \mu_2\left(J_{A_{0d}}^{1/\alpha}\right) = \frac{1}{2}\lambda_{max}\left(J_{A_{0d}}^{1/\alpha} + J^{*\,1/\alpha}_{A_{0d}}\right) \geq \frac{1}{2}\lambda_{max}\left(J_{A_{0d}} + J^{*}_{A_{0d}}\right)^{1/\alpha} = max\left(Re\,\hat{\lambda}:\hat{\lambda}\in\sigma\left(J_{A_{0d}}^{1/\alpha}\right)\right)$$

$$\geq \frac{1}{2}\lambda_{max}\left(J_{A_{0d}} + J^{*}_{A_{0d}}\right)^{1/\alpha} = max\left(Re\,\lambda^{1/\alpha}:\lambda\in\sigma(A_0)\equiv\sigma\left(J_{A_{0d}}\right)\right) \quad (3.28)$$

Thus, $A_0$ is a stability matrix if and only if $arg(\lambda)\in(-\theta_1,\theta_2)\subseteq\left(-\frac{\pi}{2},\frac{\pi}{2}\right)$; $\forall\lambda\in\sigma(A_0)$. If also $J_{A_{0d}}^{1/\alpha}$ is a stability matrix with $\mu_2\left(J_{A_{0d}}^{1/\alpha}\right)<0$, then $(1/\alpha)arg(\lambda)\in(-\theta_1/\alpha,\theta_2/\alpha)\subseteq\left(-\frac{\pi}{2},\frac{\pi}{2}\right)$ so that $|arg(\lambda)|<\frac{\alpha\pi}{2}$; $\forall\lambda\in\sigma(A_0)$, which is also a necessary condition for the fulfilment of the sufficiency-type condition (3.25) for global asymptotic stability of (2.1), which implies the stability of the matrix $J_{A_{0d}}^{1/\alpha}$ with the further constraint that $\mu_2\left(J_{A_{0d}}^{1/\alpha}\right)<0$. □

It follows after inspecting the solution (2.2), subject to (2.3)-(2.4), and Lemma 2.2 that the stability properties for arbitrary admissible initial conditions or admissible bounded controls are lost in general if $\alpha\geq 2$. However, it turns out that the boundedness of the solutions can be obtained by zeroing some of the functions of initial conditions. Note, in particular that $\varphi_j$ is required to be identically zero on its definition domain for $\overline{k-1}\cup\{0\}\ni j<\alpha-1$ $(\alpha\geq 2)$ in order that the $\Gamma$ - functions be positive (note that $\Gamma(x)$ is discontinuous at zero with an asymptote to $-\infty$ as $x\to 0^-$). This observation combined with Theorem 3.7 leads to the following direct result which is not a global stability result:

**Theorem 3.5.** Assume that $\alpha\geq 2$ and the constraint (3.25) holds with negative matrix measure $\mu_2\left(J_{A_{0d}}^{1/\alpha}\right)$. Assume also that $\varphi_j:[-h,0]\to R^n$ are any admissible functions of initial conditions for $\overline{k-1}\cup\{0\}\ni j\geq\alpha-1$ while they are identically zero if $\overline{k-1}\cup\{0\}\ni j<\alpha-1$. Then, the unforced solutions are uniformly bounded for all time independent of the delays. Also, the total solutions for admissible bounded controls are also bounded for all time independent of the delays. □

**ACKNOWLEDGMENTS**
The author is grateful to the Spanish Ministry of Education by its partial support of this work through Grant DPI2009-07197. He is also grateful to the Basque Government by its support through Grants IT378-10 and SAIOTEK S-PE09UN12.